\documentclass{amsart}     % 'we are using the compositio class'

\usepackage{amsfonts, amsthm, amssymb, amsmath, stmaryrd}
\usepackage{mathrsfs,array}
\usepackage{eucal,color,times,enumerate,accents}
\usepackage{tikz-cd}
\usepackage{url}
\usepackage{bbm}
\usepackage[utf8]{inputenc}
\usepackage{pgfplots}
\pgfplotsset{compat=newest}
\usepackage{lipsum}
\usepackage{graphicx}
\usepackage{rotating}
\usepackage{enumitem}
\usepackage{enumerate}
\usepackage{footmisc}
\usepackage{textcomp}
\usepackage[bookmarksnumbered,bookmarksopen]{hyperref}
\hypersetup{colorlinks, linkcolor=black, citecolor=black}
\usepackage{cleveref}
\usepackage{mathtools}
\usepackage{tabularx}

\newcommand{\F}{\mathbb{F}}
\newcommand{\Q}{\mathbb{Q}}
\newcommand{\Z}{\mathbb{Z}}
\newcommand{\N}{\mathbb{N}}

\newcommand{\pdiv}{\mathbf{pDiv}}

\newcommand{\mins}{\mathrm{min}}

\newcommand{\Foi}{\Fisoc^\textrm{\dag}}

\newcommand{\Qp}{\mathbb{Q}_p}

\newcommand{\et}{\textrm{\'et}}

\newcommand{\Dim}{\textrm{Dim}}

\DeclareMathAlphabet{\mathpzc}{OT1}{pzc}{m}{it}

\newcommand{\Ker}{\mathrm{Ker}}
\newcommand{\Image}{\mathrm{Im}}
\newcommand{\Coker}{\mathrm{Coker}}
\newcommand{\Ext}{\mathrm{Ext}}
\newcommand{\AJ}{\mathrm{AJ}}

\newcommand{\GL}{\mathrm{GL}}

\newcommand{\Spec}{\mathrm{Spec}}

\newcommand{\Hom}{\mathrm{Hom}}

\newcommand{\End}{\mathrm{End}}

\newcommand{\ML}{\mathrm{ML}}
\newcommand{\fl}{\mathrm{fl}}
\newcommand{\f}{\mathrm{f}}
\newcommand{\tor}{\mathrm{t}}

\newcommand{\perf}{\mathrm{perf}}

\newcommand{\crys}{\mathrm{crys}}

\newcommand{\Tr}{\mathrm{Tr}}

\newcommand{\fppf}{\mathrm{fppf}}
\newcommand{\tors}{\mathrm{tors}}
\newcommand{\tf}{\mathrm{tf}}

\newcommand{\kum}{\mathrm{Kum}}

\newcommand{\Fisoc}{\mathbf{F\textrm{-}Isoc}}
\newcommand{\SH}{\mathbf{SH}}

\newcommand*{\sExt}{\mathcal{E}\kern -.5pt xt}
\newcommand*{\sHom}{\mathcal{H}\kern -.5pt om}

\begin{document}
\newtheorem{theorem}[equation]{Theorem}
\newtheorem{hope}[equation]{Hope}
\newtheorem{proposition}[equation]{Proposition}
\newtheorem{lemma}[equation]{Lemma}
\newtheorem{claim}[equation]{Claim}
\newtheorem{corollary}[equation]{Corollary}

\newtheorem{fact}[equation]{Fact}

\theoremstyle{definition}
\newtheorem{definition}[equation]{Definition}
\newtheorem{question}[equation]{Question}
\newtheorem{conjecture}[equation]{Conjecture}
\newtheorem{answer}[equation]{Answer}
\newtheorem{remark}[equation]{Remark}
\newtheorem{example}[equation]{Example}
\newtheorem{warning}[equation]{Warning}
\newtheorem{notation}[equation]{Notation}
\newtheorem{construction}[equation]{Construction}
\title{Perfect points of abelian varieties}\author{Emiliano Ambrosi}
\email{eambrosi@unistra.fr} 
\address{Universit\'{e} de Strasbourg (IRMA)}
\thanks{The author would like to thank Anna Cadoret, Bruno Kahn, Carlo Gasbarri, Atsushi Shiho for useful discussions and Giuseppe Ancona for many suggestions on how to improve the exposition. The author is grateful to an anonymous referee whose suggestions helped in improving the exposition and the clarity of the paper and for pointing out the work of Trihan which greatly simplified the proof of Proposition \ref{overconvergencesemiabelian}. Part of this work has been done when the author was a guest of the Max Planck Institute for Mathematics in Bonn. He would like to express its gratitude to the MPIM in Bonn for its hospitality and financial support.}
\begin{abstract}
Let $p$ be a prime number, $k$ a finite field of characteristic $p>0$ and $K/k$ a finitely generated extension of fields. Let $A$ be a $K$-abelian variety such that all the isogeny factors are neither isotrivial nor of $p$-rank zero. We give a necessary and sufficient condition for the finite generation of $A(K^{\perf})$ in terms of the action of $\End(A)\otimes \Q_p$ on the $p$-divisible group $A[p^{\infty}]$ of $A$. In particular we prove that if $\End(A)\otimes \Q_p$ is a division algebra then $A(K^{\perf})$ is finitely generated. This implies the ``full'' Mordell-Lang conjecture for these abelian varieties. In addition, we prove that all the infinitely $p$-divisible elements in $A(K^{\perf})$ are torsion. These reprove and extend previous results to the non ordinary case. 
\end{abstract}

\maketitle
\section{Introduction}
Let $p$ be a prime number and $k$ a finite field of characteristic $p>0$. Let $K/k$ be a finitely generated extension of fields (e.g. $\mathbb F_p(t)/\F_p$), fix an algebraic closure $K\subseteq\overline K$ and write $K\subseteq K^{\perf}$ for the perfect closure of $K$, i.e. the smallest perfect field containing $K$ (or, equivalently, the field obtained adding to $K$ all the $p^{n}$-roots of its elements). Let $A$ be a $K$-abelian variety. Motivated by applications to the ``full'' Mordell-Lang conjecture, in this paper we study the structure of $A(K^{\perf})$ using $p$-adic cohomology. The main novelty of our approach is the use of ``mixed'' $p$-divisible groups and overconvergent F-isocrystals associated to elements in $A(K^{\perf})$. 
\subsection{Motivation}
In recent years there has been a remarkable interest in the study of the group $A(K^{\perf})$, see e.g. \cite{ambrosidaddezio, supersingular, daddeziominimal, ghiocamoosaperfection, ghiocaelliptic, rosslerperfectI, rosslerperfectII,  yuan2021}.

This interest is mainly motivated by its relation with the ``full'' Mordell-Lang conjecture (see e.g. \cite[Conjecture 1.2]{ghiocamoosaperfection}). Roughly, this conjecture states that if $\Gamma \subseteq A(\overline K)$ is a \textbf{finite rank} subgroup and $X\subseteq A_{\overline K}$ is an irreducible $\overline K$-subvariety, then $X(\overline K) \cap \Gamma$ is not Zariski dense, unless $X$ is a ``special'' (e.g. the translate of an abelian subvariety of $A$).

The characteristic zero version of the Mordell-Lang conjecture $\ML$ is a celebrated theorem of Faltings (\cite{Faltings}) for finitely generated subgroups, extended to the finite rank ones by Hindry (\cite{Hindry}). In our positive characteristic setting, the conjecture has been proved in \cite{hrushovski} under the extra assumption that $\Gamma\otimes \Z_p$ is a \textbf{finitely generated} $\Z_p$-module. However the case of arbitrary subgroups of finite rank has proven to be more elusive and few results are known.

In \cite{ghiocamoosaperfection}, Ghioca and Moosa reduced the ``full'' conjecture to the case in which the subgroup $\Gamma$ is included $A(L^{\perf})$, for $K\subseteq L$ a finite field extension. Combining this with the fact that the conjecture is known when $\Gamma$ is finitely generated, the following question arise naturally.
\numberwithin{equation}{subsection} 
\begin{question}\label{fingen}

When is $A(K^{\perf})$ finitely generated? What is the structure of $A(K^{\perf})$?
\end{question}
Our main result (Theorem \ref{main}) roughly states that whether $A(K^{\perf})$ is finitely generated or not depends only on the action of $\End(A)\otimes \Q_p$ on the $p$-divisible group of $A$ and on the $p$-rank of the isogeny factors of $A$.  As a corollary of our result, one gets the Mordell-Lang conjecture for a sufficiently generic abelian variety with Newton polygon of positive $p$-rank. To simplify the exposition, we assume for the rest of the introduction that $A$ is simple  and we refer the reader to main text (and in particular to Theorem \ref{main2}) for the general case.
\subsection{Perfect points}\label{recall}
Let us recall that, while $A(K)$ is finitely generated by the Lang-N\'eron theorem (\cite{langneron}), it is well known that $A(K^{\perf})$ is not always finitely generated. For example, if $A(K)$ contains a non-torsion element and $A$ is defined up to isogeny over $k$ or $A$ is of $p$-rank $0$, then $A(K^{\perf})$ is not finitely generated.  Even worst, Helm constructed in \cite{helmperfect} an ordinary abelian variety without isotrivial isogeny factors such that $A(K^{\perf})$ is not finitely generated. So, to have finite generation, one has to impose further conditions.

On the positive side, it is well known that the torsion subgroup $A(K^{\perf})_{\tors}\subseteq A(K^{\perf})$ is finite (see for example \cite[Page 7]{ghiocamoosaperfection}), so that the interesting part to study is its torsion free quotient $A(K^{\perf})_{\tf}:=A(K^{\perf})/A(K^{\perf})_{\tors}$. Since the $i^{th}$-power Frobenius $F^i:A\rightarrow A^{(p^i)}$ and the Verschiebung $V^i:A^{(p^i)}\rightarrow A$ induce a factorization
\begin{center}
	\begin{tikzcd}
A^{(p^i)}\arrow[bend left]{rr}{p^i}\arrow{r}{V^i}& A\arrow{r}{F^i}& A^{(p^i)} & \text{such that} & A(K^{\perf})=\bigcup_{i\in \mathbb N} A^{(p^i)}(K),
\end{tikzcd}
\end{center}
where the union is taken along the injections $F^i:A(K)\hookrightarrow A^{(p^i)}(K)$, one has that
\begin{equation}\label{1/p}
A(K)[1/p]=A(K^{\perf})[1/p].
\end{equation}
Hence, to study $A(K^{\perf})$, one is reduced to understand how much the non-torsion elements of $A(K)$ become $p^n$-divisible in $A(K^{\perf})$. There are essentially two phenomena that can make $A(K^{\perf})$ not finitely generated:
\begin{enumerate}
	\item[(a)] there might be a sequence $\{x_n\}_{n\in \N}$ of non torsion elements $x_n\in A(K)$ such that $x_n$ becomes $p^n$-divisible but not $p^{n+1}$-divisible, or
	\item[(b)] there might be a non-torsion element $x\in A(K)$ that becomes infinitely p-divisible in $A(K^{\perf})$.
\end{enumerate}
Both cases can happen and our main result says that the occurring of (a) depends only on the action of $\End(A)\otimes \Q_p$ on the $p$-divisible group of $A$ and the occurring of (b) only on the $p$-rank of $A$.
\subsection{Main results}\numberwithin{equation}{subsubsection} 
\subsubsection{Main result}
To state our main result, recall that the $p$-divisible group $A[p^{\infty}]$ of $A$ fits into a canonical connected-\'etale exact sequence
\begin{equation}\label{exactsequenceintroduction}
0\rightarrow A[p^{\infty}]^0\rightarrow A[p^{\infty}]\rightarrow A[p^{\infty}]^{\et}\rightarrow 0
\end{equation}
with $A[p^{\infty}]^0$ (resp. $A[p^{\infty}]^{\et}$) a connected (resp. \'etale) $p$-divisible group. 
Then we prove:
\begin{theorem}\label{main}
Assume that $A(K)\otimes \Q\neq 0$ (and recall that $A$ is assumed to be simple). Then:
\begin{enumerate}
\item $A(K^{\perf})$ is not finitely generated if and only if and there exists an idempotent $0\neq e\in \End(A)\otimes \Q_p$ (i.e. $e^2=e$)  that acts as $0$ on (the isogeny class of) $A[p^{\infty}]^{\et}$;
\item Every infinitely $p$-divisible point is torsion if and only if $A$ is of positive $p$-rank.
\end{enumerate}
\end{theorem}
\begin{remark}
Let us recall that, since $A$ is simple, $\End(A)\otimes \Q$ is a division algebra, hence the idempotent appearing in Theorem \ref{main}(1) has to live in $\End(A)\otimes \Q_{p}\setminus\End(A)\otimes \Q$. As often happens, it is much easier to construct $\Q_p$-linear combination of endomorphisms of $A$ (i.e. elements in $\End(A)\otimes \Q_p$) than actual endomorphisms of $A$ (i.e. elements in $\End(A)$). This kind of phenomena appears for example in the proof of the Tate conjecture for endomorphism of abelian varieties over finite fields (\cite{Tatefinitefield}). 
\end{remark}
Beyond the ordinary case, these seem to be the first general results towards the understanding of the torsion free part of $A(K^{\perf})$. Coming back to $(a)$ and $(b)$ of the previous Section \ref{recall}, Theorem \ref{main} says that case $(a)$ happens if and only if there exists an idempotent as in Theorem \ref{main}(1) and case $(b)$ happens if and only if the $p$-rank of $A$ is $0$. As an immediate corollary we get the following.
\begin{corollary}\label{maincorollary}
If $A$ has positive $p$-rank and $\End(A)\otimes \Q_p$ is a simple algebra, then $A(K^{\perf})$ is finitely generated.
\end{corollary}
Since for every Newton stratum of positive $p$-rank of the moduli space of abelian varieties of fixed dimension the generic member has $\End(A_{\overline K})\simeq \Z$, Corollary \ref{maincorollary}, together with the main results of \cite{hrushovski} and \cite{ghiocamoosaperfection}, implies the Mordell-Lang conjecture for such a generic abelian variety. 
\subsubsection{Comparison with previous results}\label{recall2}
We compare Theorem \ref{main} with some of the previously known results, assuming that ($A$ is simple and) $A(K)\otimes \Q\neq 0$. 

As already mentioned, if $A$ is isogenous to an abelian variety defined over  $k$, $A(K^{\perf})$ is not finitely generated. This is coherent with Theorem \ref{main}(1), since in this case the sequence (\ref{exactsequenceintroduction}) splits canonically up to isogeny and this splitting is induced, by the $p$-adic Tate conjecture for abelian varieties, from an idempotent $e\in \End(A)\otimes \Q_p$. Similarly, the fact that if $A$ is of $p$-rank $0$ then $A(K^{\perf})$ is not finitely generated, is coherent with Theorem \ref{main}(1), taking $e=Id_A$.

When $A$ is ordinary, Theorem \ref{main} was essentially already known, since (2) follows from \cite[Theorem 1.4]{rosslerperfectII} and (1) combining \cite[Theorem 1.1]{rosslerperfectII}) with \cite[Theorem 1.1.3]{daddeziominimal} (and their proofs). Always in the ordinary case, if $\Dim(A)\leq 2$, then $A(K^{\perf})$ is always finitely generated: this can be either deduced from \cite[Theorem 1.2 (g)]{rosslerperfectII}) or from Theorem \ref{main}(2). 
\begin{remark}
	Most of the results recalled in this section also holds replacing $k$ with $\overline k$, assuming that $A_{\overline K}$ is not isogenous to an abelian variety defined over $\overline k$. Also our Theorem \ref{main} holds replacing  $k$ with $\overline k$, as we show in Theorem \ref{corollarygeometric}, by elaborating the arguments used in the proof of Theorem \ref{main}.
	\end{remark}
\subsection{Strategy}
Our proof is mostly cohomological, in the sense that we work with $p$-divisible group and crystals. To lift our cohomological results to $\End(A)$ and  $\End(A)\otimes \Q_p$, we use the assumption that $K$ is finitely generated over a finite field, to be able to apply the $p$-adic Tate conjecture for abelian varieties.
\subsubsection{$p$-adic Abel-Jacobi maps}
To prove \Cref{main}, we start, in Section \ref{Section1}, considering various Abel-Jacobi maps. By using the short exact sequence $0\rightarrow A[p^n]\rightarrow A\xrightarrow{p^n} A\rightarrow 0$ one constructs a Abel-Jacobi map
$$\AJ:A(K)\otimes \Q\rightarrow \Ext^1(\Q_p/\Z_p,A[p^{\infty}])\otimes \Q_p.$$
Composing with the quotient map $A[p^{\infty}]\rightarrow A[p^{\infty}]^{\et}$, we get a morphism
$$\AJ^{\et}:A(K)\otimes \Q\rightarrow \Ext^1(\Q_p/\Z_p,A[p^{\infty}])\otimes \Q_p\rightarrow \Ext^1(\Q_p/\Z_p,A[p^{\infty}]^{\et})\otimes \Q_p$$
which we call the \'etale Abel-Jacobi map, and we consider its $\Q_p$-linearization 
$$\AJ^{\et}_p:A(K)\otimes \Q_{p}\rightarrow \Ext^1(\Q_p/\Z_p,A[p^{\infty}]^{\et})\otimes \Q_p,$$
which we call the p-adic \'etale Abel-Jacobi map. In Proposition \ref{keyproposition} we prove that every infinitely p-divisible element is torsion if and only if $\AJ^{\et}$ is injective and that $A(K^{\perf})$ is finitely generated if and only if $\AJ^{\et}_p$ is injective. Hence we can translate the two statements of Theorem \ref{main} into two statements on ``mixed'' $p$-divisible groups  associated to elements in $A(K)\otimes \Qp$ and $A(K)\otimes \Q$. 
\begin{remark}
Since the two properties of having  a non torsion infinitely p-divisible point and having a finitely generated group of perfect points are codified by two different maps (one $\Q_p$-linear and the other $\Q$-linear), it is natural to consider two different statements in Theorem \ref{main}. This is slightly different from what one could aspects from apparently similar motivic conjectures (see e.g. Jansen injectivity conjecture (\cite[Conj. 9.15]{jansen})). Roughly, this shows that the behavior of $\AJ^{\et}_p$ is not motivic, since  $\AJ^{\et}_p$ might not be injective even when $\AJ^{\et}$ is.
\end{remark}
\subsubsection{p-divisible groups and crystals}
For $x\in A(K)\otimes \Qp$, let 
\begin{equation}\label{fundamentalsequenceintroduction}
0\rightarrow A[p^{\infty}]\rightarrow M_x[p^{\infty}]\rightarrow \Q_p/ \Z_p\rightarrow 0 \quad \text{and}\quad 0\rightarrow A[p^{\infty}]^{\et}\rightarrow M_x[p^{\infty}]^{\et}\rightarrow \Q_p/\Z_p\rightarrow 0
\end{equation}
be the exact sequences of $p$-divisible groups representing $\AJ_p(x)$ and $\AJ^{\et}_p(x)$. By the finite generation of $A(K)$, we know that first does not split and we want to understand when and why second splits. To do this, we spread out $A\rightarrow K$ to an abelian scheme $\mathcal A\rightarrow X$ over some smooth connected $k$-variety $X$ with function field $K$ and we consider the category $\Fisoc(X)$ of F-isocrystals and the fully faithful controvariant Dieudonn\'e functor (\cite{berthelotbreenmessing2}) 
$$\mathbb D:\pdiv(X)_{\Q}\rightarrow \Fisoc(X).$$
By fully faithfulness, we translate the splitting properties of (\ref{fundamentalsequenceintroduction}) into analogous splitting properties of an exact sequence of F-isocrystals. As in \cite{ambrosidaddezio}, the advantage of doing this is that we can prove in Proposition \ref{overconvergence} that the image via $\mathbb D:\pdiv(X)_{\Q}\rightarrow \Fisoc(X),$ of the first sequence in (\ref{fundamentalsequenceintroduction}) lies inside the much better behaved subcategory $\Foi(X)\subseteq \Fisoc(X)$ of overconvergent F-isocrystals. 

Since $\mathbb D(A[p^{\infty}])$ is semisimple in $\Foi(X)$, we can apply recent advances in p-adic cohomology (\cite{tsuzukiminimal} and its improvement done in \cite{daddeziominimal}) to construct, from the splitting of $\AJ^{\et}_p$, an idempotent in $\End(A[p^{\infty}])\otimes \Q_p$ with the desired properties, which, since $K$ is finitely generated over a finite field, lifts to $\End(A)\otimes \Q_p$, by the $p$-adic  Tate conjecture for abelian varieties.

This is enough to conclude the proof of \ref{main}(1), but to complete the proof of Theorem \ref{main}(2) one needs to show that such a splitting can not exist if the sequence (\ref{fundamentalsequenceintroduction}) comes from an $x\in A(K)\otimes \Q$ and not from a random $x\in A(K)\otimes \Q_p$. This follows from  Proposition \ref{Tateproperty} which shows that even if $\End(A[p^{\infty}])$ can be big and with lots of idempotempotent, one always has that $\End(M_x[p^{\infty}])\otimes \Q_p\simeq \Q_p$ if $x\in A(K)\otimes \Q$. This is essentially due to the geometric origin of $M_x[p^{\infty}]$, which makes $M_x[p^{\infty}]$ much more rigid for a $x\in A(K)\otimes \Q$ than for a random $x\in A(K)\otimes \Q_p$. This extra rigidity is the reason for difference between the two different parts of Theorem \ref{main}.
\subsection{Organisation of the paper}
In Section \ref{Section1} we study various $p$-adic Kummer and Abel-Jacobi maps, their relation with the group of perfect points and with the extensions of $p$-divisible groups. In Section \ref{sectionproof} we use this to prove Theorem \ref{main} assuming the overconvergence result Proposition \ref{overconvergencesemiabelian}. Finally, in Section \ref{sectionoverconvergence} we prove this overconvergence result.
\vspace*{5mm}
\section{Abel-Jacobi and \'etale Abel-Jacobi maps}\label{Section1}
Let $S$ be a noetherian $\F_p$-scheme and let $A\rightarrow S$ be an abelian scheme. We write $\SH_{\fppf}(S)$ for the category of fppf sheaves in abelian groups on $S$. Write $A(S)_{\tors}\subseteq A(S)$ for the torsion subgroup of $A(S)$, $A(S)_{\tf}:=A(S)/A(S)_{\tors}$ for its torsion free quotient and 
$$A(S)_{p^{\infty}}:=\{x\in A(S) \text{ such that for every }n\in \mathbb N \text{ there exists a $y_n\in A(S)$ with $p^ny_n=x$}\}$$
for its subgroup of infinitely $p$-divisible elements.
\subsection{Kummer maps}\label{Kummersection}
\subsubsection{Kummer map}
For every $n\in \mathbb N$, the exact sequence 
$$0\rightarrow A[p^n]\rightarrow A\xrightarrow{p^n} A\rightarrow 0$$
in $\SH_{\fppf}(S)$, induces an injective morphism
$$\kum_n:A(S)/p^n\hookrightarrow H^1_{\fl}(S,A[p^n])$$
and taking the projective limit and tensoring with $\Q$, we get a commutative diagram
 \begin{center}
\begin{tikzcd}
& A(S)\otimes \Q\arrow{dl}\arrow{dr}{\kum}\arrow{d} \\
A(S)\otimes \Q_p\arrow{r}\arrow[bend right]{rr} {\kum_p}&(\varprojlim_n A(S)/p^n )\otimes \Q\arrow[hook]{r}& (\varprojlim_n H^1_{\fl}(S,A[p^n]))\otimes \Q.
\end{tikzcd}
\end{center}
We call $\kum$: $ A(S)\otimes \Q\rightarrow (\varprojlim_n H^1_{\fl}(S,A[p^n]))\otimes \Q$ the Kummer map and $\kum_p$: $A(S)\otimes \Q_p\rightarrow (\varprojlim_n H^1_{\fl}(S,A[p^n]))\otimes \Q$ the $p$-adic Kummer map. By construction, one has the following lemma, which we state for further references.
\begin{lemma}\label{basiclemma}
\begin{enumerate}
\item[]
	\item $\kum$ is injective if and only if $A(S)_{p^{\infty}}\subseteq A(S)_{\tors}$;
	\item If $A(S)_{\tf}$ is finitely generated, then $\kum_p$ is injective.
	\end{enumerate}
\end{lemma}
\proof
Statement (1) follows by tensoring with $\Q$ the short exact sequence
$$0\rightarrow A(S)_{p^{\infty}}\rightarrow A(S)\rightarrow \varprojlim_n H^1_{\fl}(S,A[p^n]).$$
For (2), one uses that if $A(S)_{\tf}$ is finitely generated, then the kernel of $A(S)\otimes \Z_p\rightarrow \varprojlim_n A(S)/p^n$
is torsion, so that the map $A(S)\otimes \Q_p\rightarrow (\varprojlim_n A(S)/p^n)\otimes \Q$ is injective. \endproof
\subsubsection{\'Etale Kummer maps}
Assume now that $S=\Spec(K)$ is the spectrum of a field and write $K^{\perf}$ for the perfection of $K$. Then, the quotient maps
$A[p^n]\rightarrow A[p^n]^{\et}$ induce a commutative diagram
\begin{center}
\begin{tikzcd}
A(K)\otimes \Q\arrow[swap]{dr}{\kum}\arrow{d}\arrow{rrd}{\kum^{\et}} \\
A(K)\otimes \Q_p\arrow[swap]{r}{\kum_p}\arrow[bend right]{rr}{\kum_p^{\et}} & (\varprojlim_n H^1_{\fl}(k,A[p^n]))\otimes \Q\arrow{r} & (\varprojlim_n H^1_{\fl}(K,A[p^n]^{\et}))\otimes \Q\simeq H^1_{\et}(K,T_p(A))\otimes \Q,
\end{tikzcd}
\end{center}
where $T_p(A):=\varprojlim_n A(\overline K)[p^n]$ is the $p$-adic \'etale module of $A$ and $H^1_{\et}(K,T_p(A))$ is its first continuous \'etale cohomology group.
We call $\kum^{\et}$: $A(K)\otimes \Q\rightarrow H^1_{\et}(K,T_p(A))\otimes \Q$ the \'etale Kummer map and $\kum^{\et}_p$: $A(K)\otimes \Q_p\rightarrow H^1_{\et}(K,T_p(A))\otimes \Q$ the $p$-adic \'etale Kummer map. The following proposition links the properties of $\kum^{\et}$ and $\kum^{\et}_p$ with the study of $A(K^{\perf})$. 
\begin{proposition}\label{keyproposition}
\begin{enumerate}
\item[]
	\item $A(K^{\perf})_{p^{\infty}}\subseteq A(K^{\perf})_{\tors}$ if and only if $\kum^{\et}$ is injective;
	\item $A(K^{\perf})_{\tf}$ is finitely generated if and only if $A(K)_{\tf}$ is finitely generated and $\kum^{\et}_p$ is injective.
		\end{enumerate}
\end{proposition}
\proof
Let us recall that 
\begin{enumerate}
	\item[(a)] Since $K\subseteq K^{\perf}$ is purely inseparable, for every finite \textbf{\'etale} group scheme $G$ the natural map 
	$H^1(K,G)\rightarrow H^1(K^{\perf},G)$ is an isomorphism (see e.g. \cite[Tag 04DZ]{stacks-project});
	\item[(b)] If $L$ is a perfect field, then 
	$H^1_{\fl}(L,H)\rightarrow H^1_{\fl}(L,H^{\et})$ is injective 
	for every finite group scheme H over $L$, since $H^1_{\fl}(L,G)=0$ for every finite \textbf{connected} group scheme $G$ (see e.g. \cite[Lemma 2.7 (a)]{kestutis}).
	\end{enumerate}
Hence $(1)$ and the only if part of $(2)$ follows from Lemma \ref{basiclemma} and the commutative diagram for $?\in \{\emptyset, p\}$:
\begin{center}
\begin{tikzcd}
A(K)\otimes \Q_?\arrow{r}\arrow{d}{\simeq}& (\varprojlim_n H^1_{\fl}(K,A[p^n]))\otimes \Q\arrow{d}\arrow{r} & H^1_{\et}(K,T_p(A)))\otimes \Q\arrow{d}{\simeq}\\
A(K^{\perf})\otimes \Q_?\arrow{r}& (\varprojlim_n H^1_{\fl}(K^{\perf},A[p^n]))\otimes \Q\arrow[hook]{r} & H^1_{\et}(K^{\perf},T_p(A)))\otimes \Q,
\end{tikzcd}
\end{center}
where the left vertical isomorphism follows from (\ref{1/p}), the right vertical isomorphism from (a) and the bottom right injection from (b).

So we are left to prove that if $A(K)_{\tf}$ is finitely generated and $\kum^{\et}_p$ is injective then $A(K^{\perf})_{\tf}$ is finitely generated. Since $A(K)_{\tf}[1/p]=A(K^{\perf})_{\tf}[1/p]$ is a finitely generated $\Z[1/p]$-module, it is enough to show that $A(K^{\perf})_{\tf}\otimes\Z_p$ is a finitely generated $\Z_p$-module. Since the kernel of $\kum_p^{\et}$ is a torsion group by assumption and $A(K)\otimes\Q=A(K^{\perf})\otimes\Q$, the group $A(K^{\perf})_{\tf}\otimes\Z_p$ injects in the torsion free quotient of the image of $\kum_p^{\et}$. Hence it is enough to show that the image of $A(K^{\perf})\otimes\Z_p$ in $H^1(K^{\perf},T_p(A))\simeq H^1(K,T_p(A))$ lies in a finitely generated sub $\Z_p$-module.

Since $A(K)_{\tf}$ is finitely generated, we can choose a set $x_1,\dots x_r\in A(K)$ which generates $A(K)_{\tf}$ and write 
$T_p(M_{x_i})$ for the $\Z_p$-linear $\pi_1(K)$-representation corresponding to the exact sequence $\kum^{\et}(x_i)$
\begin{equation}\label{equazionecheusopocodopo}
	0\rightarrow T_p(A)\rightarrow T_p(M_{x_i})\rightarrow \Z_p\rightarrow 0 \quad \text{ in } \quad H^1(K,T_p(A))\simeq \Ext^1_K(\Z_p,T_p(A))
\end{equation}
Let 
$$\Pi\subseteq \GL(T_p(M_{x_1}))\times \dots \times \GL(T_p(M_{x_r}))$$ be the image of $\pi_1(K^{\perf})$ acting on $T_p(M_{x_1})\times \dots \times T_p(M_{x_r})$ and write $K^{\perf}\subseteq L$ for the Galois extension corresponding to the closed subgroup $\Ker(\pi_1(K^{\perf})\twoheadrightarrow \Pi)$. 

Since $\Pi$ is a closed subgroup of $\GL(T_p(A))$, it is a compact $p$-adic Lie group by \cite[Corollary 9.36]{analyticprop}. In particular, by \cite[Prop. 9]{serregroupesdecongruence}, $H^1(\Pi,T_p(A))\subseteq H^1(K^{\perf},T_p(A))$ is a finitely generated $\Z_p$-module. We are left to show that the image of $A(K^{\perf})\otimes\Z_p$ in $H^1(K^{\perf},T_p(A))$ lies in $H^1(\Pi,T_p(A))$. Since $H^1(\Pi,T_p(A))$ is a sub $\Z_p$-module of $H^1(K^{\perf},T_p(A))$, it is enough to show that the image of $A(K^{\perf})$ lies in $H^1(\Pi,T_p(A))$.

The inflation-restriction exact sequence
\begin{equation}\label{diagrammafingen}
\begin{tikzcd}
0\arrow{r} & H^1(\Pi,T_p(A))\arrow{r} & H^1(\pi_1(K^{\perf}),T_p(A))\arrow{r} & H^1(\pi_1(L),T_p(A))
\end{tikzcd}
\end{equation}
reduces us to show that the composition 
$$\phi:A(K^{\perf})\rightarrow H^1(\pi_1(K^{\perf}),T_p(A))\rightarrow H^1(\pi_1(L),T_p(A))$$
 is the zero map.
Since $\pi_1(L)$ acts trivially on $T_p(M_{x_i})$, it acts trivially $T_p(A)$, so that
$$H^1(\pi_1(L),T_p(A))=\Hom(\pi_1(L),\Z_p^{p(A)})$$
is torsion free, hence it is enough to show that for every non torsion $x\in A(K^{\perf})$, there exists an $n$ such that $\phi(p^nx)=0$.
Since, by (\ref{1/p}), for every $x\in A(K^{\perf})$, there exists an $n$ such that $p^nx\in A(K)$, it is enough to show that the map 
$$\phi':A(K)_{\tf}\rightarrow H^1(\pi_1(L),T_p(A))$$
is zero.

Since $A(K)_{\tf}$ is generated by $x_1,\dots, x_r$, it is enough to show that $\phi'(x_i)=0$ for every $1\leq i\leq r$. But the exact sequence corresponding to $\phi'(x_i)$ is the restriction of the exact sequence (\ref{equazionecheusopocodopo}) to $\pi_1(L)$. By construction this sequence is an exact sequence of trivial $\pi_1(L)$-representations hence it splits as $\pi_1(L)$-module for all the $x_i\in A(K)$. Hence $\phi'(x_i)=0$ and this concludes the proof.
\endproof
\subsection{Interpretation in terms of Abel-Jacobi maps}\label{Kummersection2}
In this section, we compare the Kummer map with an Abel-Jacobi map constructed via p-divisible groups and 1-motives.

Write $\pdiv(S)$ for the category of $p$-divisible group over $S$ and $\pdiv(S)\otimes \Q$ for its isogeny category. 
\subsubsection{p-divisible group associate to a point}
Let $s\in A(S)$ be a section. Since $s:S\rightarrow A$ corresponds to a morphism of fppf $S$-groups schemes $s:\Z\rightarrow A$, we can consider the 1-motive $[s:\Z\rightarrow A]$. We now recall how to associate to $[s:\Z\rightarrow A]$ a $p$-divisible group $M_s[p^{\infty}]$ over $S$ (see for example \cite[Section 1.3]{andreattabarbieri1motivi} for more details).
Define
$$M_s[p^n]:=\frac{\Ker(s+p^n:\Z\times_S A\rightarrow A)}{\Image((p^n,-s):\Z\rightarrow \Z\times_S A)},$$ 
so that there is an exact sequence
\begin{equation}\label{fundamental exact sequence level n}
0\rightarrow A[p^n]\rightarrow M_s[p^n]\rightarrow \Z/p^n\Z_{S}\rightarrow 0
\end{equation}
of finite flat $S$-group schemes. Define 
$$M_{s}[p^{\infty}]=\varinjlim_n M_s[p^n]$$
so that $M_{s}[p^{\infty}]$ is a $p$-divisible group fitting into an exact sequence
\begin{equation}\label{fundamental exact sequence}
0\rightarrow A[p^{\infty}]\rightarrow M_s[p^{\infty}]\rightarrow (\Q_p/\Z_p)_S\rightarrow 0.
\end{equation}
We let $[M_s[p^{\infty}]]$ be the corresponding class in $\Ext^1(\Q_p/\Z_p, A[p^\infty])$.
\subsubsection{Comparison with the Kummer class}
Since $H^1_{\fl}(S,A[p^n])\simeq \Ext^1(\Z/p^{n}\Z, A[p^n])$, where the latter is the group of extension $A[p^n]$ by $\Z/p^n\Z$ as $\Z/p^n\Z$-sheaf, the Kummer map can be interpreted as a morphism
$$\kum: A(S)\rightarrow \varprojlim_n \Ext^1(\Z/p^{n}\Z, A[p^n]).$$
On the other hand, since $\Hom(\Z/p^{n}\Z, A[p^n]) $ is finite, taking $p^n$-torsion we get a natural injective morphism
$$\varphi:\Ext^1(\Q_p/\Z_p, A[p^\infty])\hookrightarrow \varprojlim_n \Ext^1(\Z/p^{n}\Z, A[p^n]).$$
In the next lemma, which follows essentially from the constructions involved, we prove that $\varphi([M_s[p^{\infty}])$ and $\kum(s)$ represent the same class
\begin{lemma}\label{comparison}
There is an equality $\kum(s)=\varphi([M_s[p^{\infty}]])$.
	\end{lemma}
\proof
It is enough to show that, for every $n$, the sequence (\ref{fundamental exact sequence level n}) identifies with the class of $\kum(s)\in H^1_{\fl}(S,A[p^n])\simeq \Ext^1(\Z/p^{n}\Z, A[p^n])$.
By definition, the $A[p^n]$-torsor $\kum(s)\in H^1_{\fl}(S,A[p^n])$ is the pullback of the inclusion of $s\hookrightarrow A$ along the multiplication by $p^n:A\rightarrow A$. 

Let $\Z/p^n\Z[\kum(s)]$ be the free $\Z/p^n\Z$-sheaf on $[\kum(s)]$, let 
$\deg:\Z/p^n\Z[\kum(s)]\rightarrow \Z/p^n\Z$ be the ``degree'' map sending $\sum n_iz_i$ to $\sum n_i$ and write $B:=\Ker(\deg)$.

By construction (see e.g. \cite[03AJ]{stacks-project}), the sequence 
\begin{equation}\label{exactsequencecheusofrapochissimo}
	0\rightarrow A[p^n]\rightarrow \widetilde{\kum(s)}\rightarrow \Z/p^n\Z\rightarrow 0,
\end{equation}
in $\Ext^1(\Z/p^{n}\Z, A[p^n])$ corresponding to $\kum(s)$, is obtained by pushing out the exact sequence
$$0\rightarrow B\rightarrow \Z/p^n\Z[\kum(s)]\rightarrow \Z/p^n\Z\rightarrow 0,$$
along the map $B\rightarrow A[p^n]$ sending the generators of the form $x-x'$ to the unique $a$ such that $x+a=x'$. 
The isomorphism of the sequence (\ref{exactsequencecheusofrapochissimo}) with the sequence (\ref{fundamental exact sequence level n}) is then induced by the map $\widetilde{\kum(s)}\rightarrow M_x[p^n]$ obtained by the universal property of pushout using the natural inclusion $A[p^n]\subseteq \{0\}\times A\subseteq \Z\times A$ and the map $\Z[\kum(s)]\rightarrow A$ sending $s\in \kum(s)$ to $(1,-s)\in \Z\times A$.
\endproof
Hence, for now on, if $S=\Spec(K)$ is the spectrum of a field, we interpret, for $?\in \{\emptyset, p\}$ and $\Delta\in \{\emptyset, \et\}$ the Kummer maps as (p-adic, \'etale) Abel-Jacobi maps
$$\AJ_?^{\Delta}: A(K)\otimes \Q_?\rightarrow \Ext^1(\Q_p/\Z_p, A[p^{\infty}]^{\Delta})\otimes \Q.$$
We can then rephrase the work done in this section in the following corollary, which is a direct consequence of Proposition \ref{keyproposition} and Lemma \ref{comparison}.
\begin{corollary}\label{corollollariocheusodavvero}
\begin{enumerate}
\item[]
	\item $A(K^{\perf})_{\tf}$ is finitely generated if and only if $A(K)_{\tf}$ is finitely generated and $\AJ^{\et}_{p}$ is injective;
 \item $A(K^{\perf})_{p^{\infty}}\subseteq A(K^{\perf})_{\tors}$ if and only if $\AJ^{\et}$ is injective.
 \end{enumerate}
\end{corollary}
\subsubsection{Rigidity of the Abel-Jacobi extension}
Suppose that $S=Spec(K)$ for a finitely generated field $K$ over $\F_p$. 
We give a first application of the interpretation of $\kum$ in terms of $p$-divisible groups, proving that the extensions in the image of $\AJ$ are more rigid than a general extension in the image of $\AJ_p$. This sets an important difference between the maps $\AJ$ and $\AJ_p$ and it is the reason why one has to consider two different statements in Theorem \ref{main}. 

We begin with an easy but important lemma, which is the only place in which some assumption on the geometry of $A$ is used.
\begin{lemma}\label{usosemplicity}
Assume that $A$ is simple and $x\in A(K)$ is a non torsion point. Then the map
$$\psi_x:\End(A)\rightarrow A(K)$$
sending $f$ to $f(x)$ is injective.
\end{lemma}
\proof
Take any morphism $f:A\rightarrow A$ such that $f(x)=0$. If $f:A\rightarrow A$ is not the zero map then, since $A$ is simple, $\Ker(f)$ is finite. On the other hand $x$ is in $\Ker(f)$ which is a contradiction with the fact that $x$ is not torsion.
\endproof
Then one has the following result, which is a consequence of the Tate conjecture for abelian varieties and a concrete incarnation of the Tate-conjecture for 1-motives.
\begin{lemma}\label{Tateproperty}
If $A$ is simple and $x\in A(K)$ is not torsion, then $\End_{\pdiv(K)}(M_{x}[p^{\infty}])\simeq \Z_p$.
\end{lemma}
\proof
Applying the functor $\Hom_{\pdiv(K)}(M_x[p^{\infty}],-)$ to the exact sequence (\ref{fundamental exact sequence}) we get an exact sequence
$$0\rightarrow \Hom_{\pdiv(K)}(M_x[p^{\infty}],A[p^{\infty}])\rightarrow \End_{\pdiv(K)}(M_x[p^{\infty}])\rightarrow \Hom_{\pdiv(K)}(M_x[p^{\infty}],\Q_p/\Z_p).$$
Since $\Hom_{\pdiv(K)}(A[p^{\infty}],\Qp/\Z_p)=0$, applying the functor $\Hom_{\pdiv(K)}(-,\Q_p/\Z_p)$ to (\ref{fundamental exact sequence}) one sees that $\Hom_{\pdiv(K)}(M_x[p^{\infty}],\Q_p/\Z_p)\simeq \End_{\pdiv(K)}(\Q_p/\Z_p)\simeq \Z_p$.
Hence, it is then enough to prove that $$\Hom_{\pdiv(K)}(M_x[p^{\infty}],A[p^{\infty}])=0.$$

Since $\Hom_{\pdiv(K)}(\Qp/\Z_p,A[p^{\infty}])=0$, applying the functor $\Hom_{\pdiv(K)}(-,A[p^{\infty}])$ to the exact sequence (\ref{fundamental exact sequence}) we get an exact sequence
$$0\rightarrow \Hom_{\pdiv(K)}(M[p^{\infty}],A[p^{\infty}])\rightarrow \End_{\pdiv(K)}(A[p^{\infty}])\rightarrow \Ext_{\pdiv(K)}^{1}(\Qp/\Z_p,A[p^{\infty}]).$$
Since $\Hom_{\pdiv(K)}(M[p^{\infty}],A[p^{\infty}])$ is torsion free, we are left to show that the natural map
$$\End_{\pdiv(K)}(A[p^{\infty}])\otimes \Q\rightarrow \Ext_{\pdiv(K)}^{1}(\Qp/\Z_p,A[p^{\infty}])\otimes \Q$$
is injective.

Consider the commutative diagram
\begin{center}
\begin{tikzcd}
\End(A)\otimes \Qp\arrow[hook]{r}{\psi_x\otimes \text{Id}}\arrow{d}{\simeq} & A(K)\otimes \Q_p\arrow{d}{\AJ_p}\\
\End_{\pdiv(K)}(A[p^{\infty}])\otimes \Q\arrow{r}& \Ext_{\pdiv(K)}^{1}(\Qp/\Z_p,A[p^{\infty}])\otimes \Q
\end{tikzcd}
\end{center}
where $\psi_x\otimes \text{Id}$ is induced by the map $\psi_x:\End(A)\rightarrow A(K)$ sending a morphism $f$ to $f(x)$. Since $K$ if finitely generated, $A(K)$ is a finitely generated group, hence by Lemma \ref{basiclemma} $\AJ_p$ is injective. By the $p$-adic Tate conjecture for abelian varieties proved in \cite[Theorem 2.6]{dejonghomomorphism}, the left vertical map is an isomorphism.
So, since $A$ is simple, we conclude by using Lemma \ref{usosemplicity}.\endproof
\subsubsection{Kummer class and semiabelian schemes}\label{duality1motives}
Let $s\in A(S)$. As a second application of the interpretation of $\kum$ in terms of $p$-divisible groups we give a geometric interpretation of the  Cartier dual of the class of $[M_x[p^{\infty}]]$. This will be important to prove Proposition \ref{overconvergence}. The dual of the 1-motive $[\Z\rightarrow A]$ is a semiabelian scheme $$0\rightarrow \mathbb G_{m,S}\rightarrow G_s\rightarrow  A^{\vee}\rightarrow 0,$$
where $A^{\vee}$ is the dual abelian variety, and the $p$-divisible group $ G_s[p^{\infty}]$ of $G_s$ is the Cartier dual $M_{s}[p^{\infty}]^{\vee}$ of $M_{s}[p^{\infty}]$ (see for example \cite[Section 1.3]{andreattabarbieri1motivi}).
Hence, the class  of the dual of the extension (\ref{fundamental exact sequence}) in $\Ext^1(A^{\vee}[p^{\infty}],\mu_{p^{\infty}})$ is the extension
$$0\rightarrow \mathbb G_{m}[p^{\infty}]\rightarrow G_s[p^{\infty}]\rightarrow A^{\vee}[p^{\infty}]\rightarrow 0,$$
associated to the $p$-divisible group of a semi-abelian $S$-scheme $G_s\rightarrow S$.
\section{On the injectivity of the \'etale Abel-Jacobi map}\label{sectionproof}
\numberwithin{equation}{subsection} 
In this section we prove the main theorem of the paper (Theorem \ref{main2}) and its geometric variant (Theorem \ref{corollarygeometric}) assuming an overconvergence result (Proposition \ref{overconvergencesemiabelian}) which will be proved in the next Section \ref{sectionoverconvergence} (since it relies on different techniques). 
\subsection{Notation and statements}\label{notation and statements}
We assume that $k$ is a finite field, $K/k$ is a finitely generated field extension and $A$ a $K$-abelian variety. Write $p(A)$ (resp. $r(A)$) for the $p$-rank of $A$ (resp. the rank of $A(K)$, which is finite by the Lang-N\'eron theorem) and if $A_1,\dots, A_n$ are the simple isogeny factors of $A$, set $p(A)^{\mins}$ (resp. $r(A)^{\mins}$) as the minimum of $p(A_i)$ (resp. of $r(A_i)$). If $e\in \End(A)\otimes \Qp$, we write $e[p^{\infty}]\in \End(A[p^{\infty}])\otimes \Q_p$ (resp. $e[p^{\infty}]^{\et}\in \End(A[p^{\infty}]^{\et})\otimes \Q$) for the induced morphism. Finally, set 
$$A(S)_{p^{\infty}}:=\{x\in A(S) \text{ such that for every }n\in \mathbb N \text{ there exists a $y_n\in A(S)$ with $p^ny_n=x$}\}.$$
In this section we prove the following.
\begin{theorem}\label{main2}
Assume that $r(A)^{\mins}>0$. Then:
\begin{enumerate}
\item $A(K^{\perf})$ is not finitely generated if and only if and there exists an idempotent $0\neq e\in \End(A)\otimes \Q_p$ (i.e. $e^2=e$) such that $0=e[p^{\infty}]^{\et}\in \End(A[p^{\infty}]^{\et})\otimes \Q_p$;
\item $A(K^{\perf})_{p^{\infty}}\subseteq A(K^{\perf})_{\tors}$ if and only if $p(A)^{\mins}>0$ .
\end{enumerate}
\end{theorem}
Since $A(K^{\perf})_{\tors}$ is finite by \cite[Page 7]{ghiocamoosaperfection}, thanks to Corollary \ref{corollollariocheusodavvero}, Theorem \ref{main} is equivalent to the following.
\begin{theorem}\label{intro injectivity etale Abel-jacobi}
Assume that $r(A)^{\mins}>0$. Then:
\begin{enumerate}
\item The morphism 
$$\AJ^{\et}_p: A(K)\otimes \Q_p\rightarrow \Ext^1_{\pdiv(K)}(\Qp/\Z_p, A[p^{\infty}]^{\et})\otimes\Q$$
is not injective if and only there exists an idempotent $0\neq e\in \End(A)\otimes \Q_p$ such that $0=e[p^{\infty}]^{\et}\in \End(A[p^{\infty}]^{\et})\otimes \Q$;
\item The morphism
$$\AJ^{\et}: A(K)\otimes \Q\rightarrow \Ext^1_{\pdiv(K)}(\Qp/\Z_p, A[p^{\infty}]^{\et})\otimes \Q$$
is not injective if and only if $p(A)^{\mins}=0$.
\end{enumerate}
\end{theorem}
\numberwithin{equation}{subsubsection} 
\subsection{Preliminaries and the first implication}\label{preliminaries}
\subsubsection{Reduction to A simple}
Since the assumptions and the conclusions are stable by products and isogenies of abelian varieties, we can assume that $A$ is simple (that will be used to apply Lemmas \ref{usosemplicity} and \ref{Tateproperty}) and $r(A)>0$. Since the statements with $p(A)=0$ are trivial, we can assume that $p(A)>0$. 
\subsubsection{First implication}
We first prove the if part of Theorem \ref{intro injectivity etale Abel-jacobi}(1). Assume that there exists an idempotent $0\neq e\in \End(A)\otimes \Q_p$ such that $e[p^{\infty}]^{\et}=0$ in $\End(A[p^{\infty}]^{\et})\otimes \Q_p$ . Chose an $n$ such that $p^ne=:u\in \End(A)\otimes \Z_p$. Since $e[p^{\infty}]^{\et}=0$ and $\End(A[p^{\infty}]^{\et})$ is torsion free, also $u[p^{\infty}]^{\et}=0$. Take a non torsion $x\in A(K)$ (which exists by assumption).
Since $A$ is simple, by Lemma \ref{usosemplicity}, the map
$$\psi_x\otimes \text{Id}_{\Q_p}:\End(A)\otimes \Q_p\rightarrow A(K)\otimes \Q_p$$
is injective, where $\psi_x:\End(A)\rightarrow A(K)$ is the map sending $f$ to $f(x)$. Hence $e(x)\neq 0$ therefore $u(x)\neq 0$. 
The commutative diagram
\begin{center}
\begin{tikzcd}
A(K)\otimes \Qp\arrow{r}{\AJ_p^{\et}}\arrow{d}{u} & \Ext_{\pdiv(K)}^{1}(\Qp/\Z_p,A[p^{\infty}]^{\et})\otimes \Q\arrow{d}{u=0}\\
A(K)\otimes \Qp\arrow{r}{\AJ_p^{\et}}\arrow{r}& \Ext_{\pdiv(K)}^{1}(\Qp/\Z_p,A[p^{\infty}]^{\et})\otimes \Q
\end{tikzcd}
\end{center}
shows that $u(x)$ goes to zero in $\Ext_{\pdiv(K)}^{1}(\Qp/\Z_p,A[p^{\infty}]^{\et})\otimes \Q$ . This concludes the proof of the if part of Theorem \ref{intro injectivity etale Abel-jacobi}(1).
\subsubsection{Reduction to Proposition \texorpdfstring{\ref{keylemma}}-}
We are left to prove the only if part of Theorem \ref{intro injectivity etale Abel-jacobi}(1) and \ref{intro injectivity etale Abel-jacobi}(2). We first show that the following Proposition \ref{keylemma} implies Theorem \ref{intro injectivity etale Abel-jacobi}. 
\begin{proposition}\label{keylemma}
Let $x\in A(K)\otimes \Q_p$ be such that $\AJ_p(x)=0$. Then there exists an idempotent $0\neq e\in \End(M_x[p^{\infty}])\otimes \Qp$ which preserves the sub $p$-divisible group  $A[p^{\infty}]\subseteq M_x[p^{\infty}]$ and it induces a non-zero idempotent $e[p^{\infty}]\in \End(A[p^{\infty}])\otimes \Qp$ acting as $0$ on $A_x[p^{\infty}]^{\et}$.
\end{proposition}
Assume that Proposition \ref{keylemma} holds. Then Theorem \ref{intro injectivity etale Abel-jacobi}(2) follows from it and Lemma \ref{usosemplicity}. To deduce Theorem \ref{intro injectivity etale Abel-jacobi}(1), we use that, by the $p$-adic Tate conjecture for abelian varieties proved in \cite[Theorem 2.6]{dejonghomomorphism}, the natural map
$$\End(A)\otimes \Q_p \xrightarrow{\simeq} \End_{\pdiv(K)}(A[p^{\infty}])\otimes \Q_p$$
an isomorphism, so that $e[p^{\infty}]$ is induced by a non-zero idempotent in $\End(A)\otimes \Q_p$ acting as $0$ on $A_x[p^{\infty}]^{\et}$. Hence we are left to prove Proposition \ref{keylemma}.
\subsection{Proof of Proposition \texorpdfstring{\ref{keylemma}}-}
\subsubsection{Spreading out}\label{mainabel2}
Let $x\in A(K)\otimes \Q_p$ be such that $\AJ_p(x)=0$. To prove Proposition \ref{keylemma} we can replace $x$ with $p^nx$ hence we may and do assume that $x\in A(K)\otimes \Z_p$ is not torsion.
Let
\begin{equation}\label{exact sequence che uso dappertutto}
0\rightarrow A[p^{\infty}]\rightarrow M_{x}[p^{\infty}]\rightarrow \Q_p/\Z_p\rightarrow 0
\end{equation}
and 
\begin{equation}\label{exact sequence che uso dappertutto etale}
0\rightarrow A[p^{\infty}]^{\et}\rightarrow M_{x}[p^{\infty}]^{\et}\rightarrow \Q_p/\Z_p\rightarrow 0
\end{equation}
be the extensions associated to $\AJ_p(x)$ and $AJ^{\et}_p(x)$ respectively. Since $AJ^{\et}_p(x)=0$, the exact sequence (\ref{exact sequence che uso dappertutto etale}) splits. Replacing $k$ with a finite field extension, we can assume that $k$ is algebraically closed in $K$ and take an affine smooth geometrically connected $k$-variety $X$ with function field $K$. Replacing $X$ with a dense open subset, we can assume that $A$ extends to an abelian scheme $\mathcal A\rightarrow X$ with constant Newton polygon and that, since $A(K)$ is finitely generated, the natural map $\mathcal A(X)\otimes \Z_p\rightarrow A(K)\otimes \Z_p$ is an isomorphism. In particular, $x$ extends to a non-torsion element $\mathfrak t\in \mathcal A(X)\otimes \Z_p$.
By \cite{dejongcrystallineviaformal}, the natural functor 
$$\pdiv(X)\otimes \Q\rightarrow \pdiv(K)\otimes \Q$$
is fully faithful, so that our assumption is equivalent to the fact that the sequence 
\begin{equation}\label{sequence spreading out}
0\rightarrow \mathcal A[p^{\infty}]^{\et}_{X}\rightarrow \mathcal M_{\mathfrak t}[p^{\infty}]^{\et}_{X}\rightarrow \Q_p/\Z_p\rightarrow 0
\end{equation}
splits in $\pdiv(X)\otimes \Q$ and we know (by Lemma \ref{basiclemma}) that the exact sequence
\begin{equation}\label{sequence spreading outnonsplit}
0\rightarrow \mathcal A[p^{\infty}]_{X}\rightarrow \mathcal M_{\mathfrak t}[p^{\infty}]_{X}\rightarrow \Q_p/\Z_p\rightarrow 0
\end{equation}
does not split.
\subsubsection{F-isocrystals}\label{mainabel3}
Let $\Fisoc(X)$ be the category of F-isocrystals over $X$ (as defined for example in \cite[Section A.1]{morrowvariational}). By \cite[Corollary 4.2]{kedlayanotes}, every F-isocrystals $\mathcal E$ with constant Newton polygon admits a slope filtration 
$$0=\mathcal E_s\subseteq \mathcal E_{s+1}\dots \subseteq \mathcal E_{r-1}\subseteq \mathcal E_{r}=\mathcal E$$
such that $\mathcal E_i/\mathcal E_{i-1}$ is isoclinic of some slope $s_i\in \Q$ with $s_i<s_{i+1}$.
By \cite{berthelotbreenmessing2}, there is a fully faithful controvariant functor $\mathbb D: \pdiv(X)\otimes \Q\rightarrow \Fisoc(X)$. Write  
$$\mathcal E:=\mathbb D(\mathcal A[p^{\infty}]);\quad  \mathcal O^{\crys}_X:=\mathbb D(\mathbb \Q_p/\Z_p); \quad \mathcal E_{\mathfrak t}:=\mathbb D(\mathcal M_{\mathfrak t}[p^{\infty}]_{X}),$$
so that $\mathcal E$ and $\mathcal E_{\mathfrak t}$ have constant Newton polygon by the preliminary reduction. Recall that the slopes appearing in a F-isocrystal associated to a $p$-divisible group are between $0$ and $1$ and that a $p$-divisible group is \'etale if and only after applying $\mathbb D$ has constant slope $0$. Hence
$$\mathcal E_{1}:=\mathbb D(\mathcal A[p^{\infty}]^{\et}_X) \quad \text{and}\quad \mathcal E_{\mathfrak t,1}:=\mathbb D(\mathcal M_{\mathfrak t}[p^{\infty}]^{\et}_X) $$ are the sub $F$-isocrystals of minimal slope of $\mathcal E$ and $\mathcal E_{\mathfrak t}$, respectively. 
Then the sequences (\ref{sequence spreading out}) and (\ref{sequence spreading outnonsplit}) are sent to exact sequences

\begin{equation}\label{exact sequence crystals non-split}
0\rightarrow \mathcal O^{\crys}_{X}\rightarrow \mathcal E_t \rightarrow \mathcal E\rightarrow 0 \quad \text{and}
\end{equation}
\begin{equation}\label{exact sequence crystals}
0\rightarrow \mathcal O^{\crys}_{X}\rightarrow \mathcal E_{t,1} \rightarrow \mathcal E_1\rightarrow 0.
\end{equation}
By fully faithfulness of $\mathbb D:\pdiv(X)\otimes \Q\rightarrow \Fisoc(X)$ and the assumption, the sequence (\ref{exact sequence crystals}) splits and (\ref{exact sequence crystals non-split}) does not split.
\subsubsection{Overconvergence}\label{mainabel4}
Let $\Foi(X)$ be the category of overconvergent F-isocrystals over $X$ (see for example \cite[Definition 2.3.6]{berthelotrigide}). By \cite[Theorem 2.4.2]{berthelotrigide}, every $F$-isocrystals is convergent, hence there is a natural functor $\Phi:\Foi(X)\rightarrow \Fisoc(X)$.

Recall that, by \cite{kedlayafully}, the functor $\Phi:\Foi(X)\rightarrow \Fisoc(X)$ is fully faithful, so that we can identify $\Foi(X)$ with a full subcategory of $\Fisoc(X)$. If $\mathcal G$ in $\Fisoc(X)$ is in the essential image of $\Phi:\Foi(X)\rightarrow \Fisoc(X)$ we say that it is overconvergent and we write $\mathcal G^{\dagger}$ for its (unique) overconvergent extension.
By \cite{etesseschemiab}, $\mathcal E$ is overconvergent. As a consequence of Proposition \ref{overconvergencesemiabelian}, that will be proved in Section \ref{sectionoverconvergence}, and the geometric interpretation of $\mathcal E_{\mathfrak t}$ given in Section \ref{duality1motives}, we can show that $\mathcal E_{\mathfrak t}$ is also overconvergent.
\begin{proposition}\label{overconvergence}
The F-isocrystal $\mathcal E_{\mathfrak t}$ is overconvergent.
\end{proposition}
\proof
Since inside $\Ext^1_{\Fisoc(X)}(\mathcal E, \mathcal O^{\crys}_X)$ the class of $\mathcal E_{\mathfrak t}$ is a $\Q_p$-linear combination of classes $\mathcal E_{\mathfrak v}$ with $\mathfrak v\in \mathcal A(X)$ and the morphism $\Ext^1_{\Foi(X)}(\mathcal E^{\dagger}, \mathcal O^{\dagger}_X)\rightarrow \Ext^1_{\Fisoc(X)}(\mathcal E, \mathcal O^{\crys}_X)$ is $\Q_p$-linear, we can assume that $\mathfrak t\in \mathcal A(X)$.
It is then enough to show that $\mathcal E^{\vee}_{\mathfrak t}(1)$ (where $(-)^{\vee}$ is the dual F-isocrystals and $(-)(1)$ is the Tate twist) is overconvergent. By Section \ref{duality1motives} and the compatibility of the functor $\mathbb D$ with dualities  (\cite[(5.3.3.1)]{berthelotbreenmessing2}), one has that $\mathcal E^{\vee}_{\mathfrak t}(1)$ identifies with $\mathbb D(G[p^{\infty}])$, where $G[p^{\infty}]$ is the p-divisible group of an algebraic group $G$ which is an extension 
$$0\rightarrow \mathbb G_m\rightarrow G\rightarrow A\rightarrow 0$$
of an abelian variety and a $\mathbb G_m$. Then the overconvergence of  $\mathcal E^{\vee}_{\mathfrak t}(1)$ follows from Proposition \ref{overconvergencesemiabelian}, that we will prove in the next Section \ref{sectionoverconvergence}.
\endproof
Since $\mathcal E$ and $\mathcal E_{\mathfrak t}$ are overconvergent and the functor $\Phi:\Foi(X)\rightarrow \Fisoc(X)$ is fully faithful, the non-split exact sequence (\ref{exact sequence crystals non-split}) lifts to a non-split exact sequence 
\begin{equation}\label{overconvergentexact}
0\rightarrow \mathcal O^{\dagger}_X\rightarrow \mathcal E^{\dagger}_{\mathfrak t}\xrightarrow{\pi} \mathcal E^{\dagger}\rightarrow 0.
\end{equation}
On the other hand, by construction, the exact sequence (\ref{exact sequence crystals}) is obtained by applying $\Phi:\Foi(X)\rightarrow \Fisoc(X)$ to (\ref{overconvergentexact}) and then base changing it along $\mathcal E_1\rightarrow \mathcal E$.

%, so that the sequence (\ref{exact sequence crystals}) makes the following diagram with exact rows cartesian:
%\begin{center}
%\begin{tikzcd}
%0\arrow{r}&\mathcal O^{\crys}_X \arrow{r}\arrow[equal]{d} &\mathcal E_{t,1}\arrow{r}\arrow{d}\arrow[phantom]{dr}{\Box}& \mathcal E_1\arrow{r}\arrow[swap]{d} &0 \\
%
%0\arrow{r}&\mathcal O^{\crys}_X \arrow{r} &\mathcal E_t \arrow{r}& \mathcal E\arrow{r}& 0,
%\end{tikzcd}
%\end{center}
%where the vertical maps are the natural inclusions. 
\subsubsection{Minimal slope conjecture}\label{mainabel5}
Chose a splitting $s:\mathcal E_1\rightarrow \mathcal E_{\mathfrak t,1}$ of the sequence (\ref{exact sequence crystals}). Consider the smallest overconvergent object $\widetilde{\mathcal E}^{\dagger}$ contained in $\mathcal E^{\dagger}_{\mathfrak t}$ and containing $s(\mathcal E_1)$. 

Since $p(A)>0$, we have $\widetilde{\mathcal E}^{\dagger}\neq 0$. By the recent work \cite{tsuzukiminimal} and its improvement done in \cite[Theorem 4.1.3]{daddeziominimal}, one has $s(\mathcal E_1)=\widetilde{\mathcal E}_1$
so that $\widetilde{\mathcal E}^{\dagger}\cap \mathcal O^{\crys}_{X}=0$. Hence the natural composite map 
$$\widetilde{\mathcal E}^{\dagger}\hookrightarrow \mathcal E^{\dagger}_{\mathfrak t}\xrightarrow{\pi} \mathcal E^{\dagger}$$
is injective and it induces an isomorphism $\pi: \widetilde{\mathcal E}^{\dagger}\xrightarrow{\simeq} \pi(\widetilde{\mathcal E}^{\dagger})$.
By construction, the sequence (\ref{overconvergentexact}) splits  after base change along $\pi(\widetilde{\mathcal E}^{\dagger})\subseteq \mathcal E^{\dagger}$.

Since the sequence (\ref{overconvergentexact}) does not split, $\pi(\widetilde{\mathcal E}^{\dagger})\neq \mathcal E^{\dagger}$. By a result of P\'al (\cite[Theorem 1.2]{Pal}), the overconvergent F-isocrystals $\mathcal E^{\dagger}$ is semisimple hence there is a projection $\widetilde e:\mathcal E^{\dagger}\rightarrow \mathcal E^{\dagger}$ onto $\pi(\widetilde{\mathcal E}^{\dagger})$. Since $\widetilde{\mathcal E}$ contains $s(\mathcal E_1)$, the non-zero idempotent $1-\widetilde e$ acts as zero on $\mathcal E_1$. By the faithfulness of the composite functor $\Foi(X)\xrightarrow{\Phi} \Fisoc(X)\xrightarrow{\mathbb D} \pdiv(X)\otimes \Q$, we get a a non-zero idempotent $e[p^{\infty}]$ in $\End(A[p^{\infty}])\otimes \Qp$ acting as zero on $A[p^{\infty}]^{\et}$. Observe that the composite map 
$$\mathcal E_{\mathfrak t}^{\dagger}\xrightarrow{\pi} \mathcal E^{\dagger}\xrightarrow{\widetilde e}\pi(\widetilde{\mathcal E}^{\dagger})\xrightarrow{\pi^{-1}}\widetilde{\mathcal E}^{\dagger}\subseteq \mathcal E_{\mathfrak t}^{\dagger}$$
is a projection onto $\widetilde{\mathcal E}^{\dagger}$. Hence there exists a  non-zero idempotent $e\in \End(\mathcal E_{\mathfrak t}^{\dagger})\simeq \End_{\pdiv(K)}(M_x[p^{\infty}])\otimes \Qp$
 which induces the non-zero idempotent in $e[p^{\infty}]\in \End(A[p^{\infty}])\otimes \Qp$ acting as $0$ on $A_x[p^{\infty}]^{\et}$. This concludes the proof of Proposition \ref{keylemma}.
\subsection{Geometric variant}\label{sectionproof3}
\numberwithin{equation}{subsection} 
Write $L:=\overline kK\subseteq \overline K$ for the field generated by $\overline k$ and $K$ in $\overline K$. Let $\Tr_{\overline K/\overline k}(A)$ be the $(\overline K/\overline k)$-trace of $A_{\overline K}$ (i.e. the biggest $\overline k$-isotrivial quotient $A_{\overline K}\rightarrow\Tr_{\overline K/\overline k}(A)$ of $A_{\overline K}$). A modification of the previous arguments gives us the following geometric variant.
\begin{theorem}\label{corollarygeometric}
Assume that $r(A)^{\mins}>0$. Then:
	\begin{enumerate}
	\item If $\Tr_{\overline K/\overline k}(A)=0$, then $A(L^{\perf})$ is not finitely generated if and only if there exists an idempotent $0\neq e\in \End(A_{L})\otimes \Q_p$ such that $0=e[p^{\infty}]^{\et}\in \End(A_L[p^{\infty}]^{\et})\otimes \Q_p$;
	\item $A(L^{\perf})_{p^{\infty}}\subseteq A(L^{\perf})_{\tors}$ if and only if $p(A)^{\mins}>0$.
\end{enumerate}
\end{theorem}
\proof
Since $A(L)_{\tf}$ is finitely generated by the Lang-N\'eron theorem and the action of $\pi_1(K)$ on $\End(A)$ factors through a finite quotient, there exists a finite extension $K\subseteq K'\subseteq L$ such that $A(K')\otimes \Q=A(L)\otimes \Q$ and $\End(A_{K'})=\End(A_L)$. For $?\in \{\emptyset, p\}$ we consider the commutative diagrams
 \begin{center}
\begin{tikzcd}
A(K')\otimes \Q_? \arrow{r}\arrow{d}{\simeq} & \Ext^1_{\pdiv(K')}(\Q_p/\Z_p,A[p^{\infty}]^{\et})\otimes \Q \arrow{d}\arrow{r}{\simeq}&H^1(\pi_1(K'),T_p(A))\otimes \Q\arrow{d}{\phi} \\
A(L)\otimes \Q_?\arrow{r} & \Ext^1_{\pdiv(L)}(\Q_p/\Z_p,A[p^{\infty}]^{\et})\otimes \Q \arrow{r}{\simeq}& H^1(\pi_1(L),T_p(A))\otimes \Q
\end{tikzcd}
\end{center}
Moreover, if $\Tr(A)=0$ then $A(L^{\perf})_{\tors}$ is finite by \cite{ambrosidaddezio}. Since $A(L)_{\tf}$ is finitely generated, by Corollary \ref{corollollariocheusodavvero} and Theorem \ref{intro injectivity etale Abel-jacobi} it is enough to show that $\phi$ is injective. 
Since $\pi_1(L)\subseteq \pi_1(K')$ is an normal subgroup, the Hochschild-Serre spectral sequence gives us an exact sequence
$$0\rightarrow H^1(\pi_1(K')/\pi_1(L),T_p(A)^{\pi_1(L)})\otimes \Q\rightarrow H^1(\pi_1(K'),T_p(A))\otimes \Q \rightarrow H^1(\pi_1(L),T_p(A))\otimes \Q.$$
Since $\pi_1(K')/\pi_1(L)$ is pro-cyclic , one has
$$H^1(\pi_1(K')/\pi_1(L),T_p(A)^{\pi_1(L)})\otimes \Q \simeq 
(T_p(A)^{\pi_1(L)}\otimes \Q)_{\pi_1(K')/\pi_1(L)}$$
where the last term are the coinvariants. But since $A(K^{\perf})[p^{\infty}]$ is finite, one has
$$(T_p(A)^{\pi_1(L)}\otimes \Q)^{\pi_1(K')/\pi_1(L)}=(T_p(A)\otimes \Q)^{\pi_1(K')}=0=(T_p(A)^{\pi_1(L)}\otimes \Q)_{\pi_1(K')/\pi_1(L)},$$
and this concludes the proof. \endproof
\section{Overconvergence}\label{sectionoverconvergence}
\numberwithin{equation}{subsection} 
\subsection{Statement}
Let $X$ be a smooth geometrically connected variety over a finite field $k$ of characteristic $p$ and let \begin{equation}\label{definingsemiabelian}
0\rightarrow W\rightarrow G\rightarrow A\rightarrow 0
\end{equation}
be an extension of an abelian $X$-scheme $A$ by a torus $W$ over $X$.
By applying the Dieudonn\'e functor $\mathbb D:\pdiv(X)\otimes \Q\rightarrow \Fisoc(X)$ to the exact sequence
$0\rightarrow W[p^{\infty}]\rightarrow G[p^{\infty}]\rightarrow A[p^{\infty}]\rightarrow 0$,
we get an exact sequence
$$0\rightarrow \mathbb D(A[p^{\infty}])\rightarrow \mathbb D(G[p^{\infty}])\rightarrow \mathbb D(W[p^{\infty}])\rightarrow 0.$$
The main result of this section is the following.
\begin{proposition}
\label{overconvergencesemiabelian}
	The $F$-isocrystal $\mathbb D(G[p^{\infty}])$ is overconvergent.
\end{proposition}
To prove Proposition \ref{overconvergencesemiabelian}, we reduce to the case in which $X$ is a curve and the abelian scheme has everywhere semistable reduction. Then, in Section \ref{pdiv}, we use a result of Trihan (\cite{Trihan}) to reduce to prove a semistability result for $G[p^{\infty}]$. We conclude the proof in Sections \ref{construction} and \ref{conclusion}, proving this semistability.
\subsection{Preliminary reductions}\label{preliminary reduction}
By \cite[Lemma 4.2]{cutbycurve}, to prove overconvergence, we can freely replace $X$ with a smooth variety $Y$ admitting a dominant morphism $Y\rightarrow X$. So we can assume that $W\simeq \mathbb G_{m,X}^m$ and
that $A(X)[n]\simeq (\mathbb Z/n\Z)^{2g}$ for some fixed $n\geq 3$ coprime with $p$. By \cite[Proposition 4.7, Expos\'e IX, Pag. 48]{SGA7-II}, this last condition implies that, for every smooth curve $C$ and every morphism $C\rightarrow X$, the abelian scheme $A\times_X C$ has everywhere semistable reduction. Moreover, by de Jong's alteration theorem (\cite{dejonghomomorphism}), we can assume that $X$ admits a compactification whose complementary is a normal crossing divisor. In this situation, by \cite[2.5]{dejonghomomorphism} and \cite[Corollary 3.14]{Trihan}, for every smooth curve $C$ and every morphism $f:C\rightarrow X$, the F-isocrystals $f^*\mathbb D(A[p^{\infty}])\simeq \mathbb D(A\times_XC[p^{\infty}])$ has everywhere semistable reduction. Therefore we can apply the cut by curve criterion for overconvergence proved in \cite[Lemma 6.7]{cutbycurve} to reduce to the case in which $X$ is a curve. So from now we assume that $X$ is a curve with smooth compactification $\overline X$ and $A$ has every everywhere semistable reduction.
\subsection{Passing to p-divisible groups}\label{pdiv}
For every $x\in \overline X-X$ we let $S_x$ be the spectrum of the completion of $\overline X$ in $x$ and $\eta_x$ the generic point of $S_x$. Write  $A_{\eta_x}$ and  $G_{\eta_x}$ for the base change of $A\rightarrow X$ and $G  \rightarrow  X$ trough $\eta_x \rightarrow X$. By \cite[Theorem 4.5]{Trihan}, to prove Proposition \ref{overconvergencesemiabelian}, it is enough to show that for every $x\in \overline X-X$, the p-divisible group $G_{\eta_x}[p^{\infty}]$ is semistable, i.e. that there exists a filtration 
$$G_{\eta_x}[p^{\infty}]^{\tor}\subseteq G_{\eta_x}[p^{\infty}]^{\f}\subseteq G_{\eta_x}[p^{\infty}]$$
such that:
\begin{enumerate}
\item $G_{\eta_x}[p^{\infty}]^{\f}$ and $G_{\eta_x}[p^{\infty}]/G_{\eta_x}[p^{\infty}]^{\tor}$ extend to p-divisible groups $G[p^{\infty}]_{x,1}$ and $G[p^{\infty}]_{x,2}$ over $S_x$. In this case, by \cite{dejonghomomorphism}, the natural map $G_{\eta_x}[p^{\infty}]^{\f}\rightarrow G_{\eta_x}[p^{\infty}]/G_{\eta_x}[p^{\infty}]^{\tor}$ extends to a map $G[p^{\infty}]_{x,1}\rightarrow G[p^{\infty}]_{x,2}$;
\item $\Ker(G[p^{\infty}]_{x,1}\rightarrow G[p^{\infty}]_{x,2})$ is a multiplicative p-divisible group and $\Coker(G[p^{\infty}]_{x,1}\rightarrow G[p^{\infty}]_{x,2})$ is an \'etale p-divisible group.
\end{enumerate}
Since the situation is now entirely local, we drop the subscript $x$ from the notation.
\subsection{Construction of the filtration}\label{construction}
By \cite[Proposition 7, Pag. 292]{neronmodels} and its proof, there exists an exact sequence of smooth group $S$-schemes with connected fibers
$$ 0\rightarrow W\rightarrow \mathcal G^0\rightarrow \mathcal A^0\rightarrow 0,$$
where $\mathcal A\rightarrow S$ be the N\'eron-model of $A_\eta$ and $\mathcal A^0\rightarrow S$ is its connected component of the identity, having as generic fiber the sequence 
$$0\rightarrow W_\eta\rightarrow G_\eta\rightarrow A_\eta\rightarrow 0.$$
Since $A_{\eta}$ has semistable reduction, the special fiber $\mathcal A^0_s$ fits into an exact sequence
$$0\rightarrow T\rightarrow \mathcal A_s^0\rightarrow B\rightarrow 0$$
with $T$ a $k$-torus and $B$ a $k$-abelian variety. 
Let $\mathcal A^0[p^n]^{\f}\subseteq \mathcal A^0[p^n]$ be the maximal subgroup which is finite over $S$ and $\mathcal A^0[p^n]^{\tor}\subseteq \mathcal A^0[p^n]^{\f}$ be the unique lifting of the finite subgroup $T[p^n]\subseteq \mathcal A^0_{s}[p^n]$ to $\mathcal A^0[p^n]^{\f}$. For  $?\in\{\tor,\f\}$, we define  $\mathcal G^0[p^n]^{?}\subseteq \mathcal G^0[p^n]$ via the following cartesian diagram with exact rows:
\begin{equation}\label{diagramma}
\begin{tikzcd}
0\arrow{r} & W[p^n] \arrow{r}\arrow[equal]{d} & \mathcal G^0[p^n]^{?}\arrow{r}\arrow{d}\arrow[phantom]{rd}{\Box}  & \mathcal A^0[p^n]^{?}\arrow{r}\arrow{d} & 0\\
0\arrow{r} & W[p^n] \arrow{r} & \mathcal G^0[p^n]\arrow{r}  & \mathcal A^0[p^n]\arrow{r}&  0.
\end{tikzcd}
\end{equation}
Taking the direct limit with $n$ and applying \cite[Proposition 5.6, Expos\'e IX, Pag. 180]{SGA7-II} and \cite[(2.4.3)]{messingcrystals}, we get a filtration $\mathcal G^0[p^{\infty}]^{\tor}\subseteq \mathcal G^0[p^{\infty}]^{\f}\subseteq \mathcal  G^0[p^{\infty}],$ of $p$-divisible groups.
Set 
$$G[p^{\infty}]_{\eta}^{\tor}:=\mathcal G^0[p^{\infty}]^{\tor}_{\eta};\quad G[p^{\infty}]_{\eta}^{\f}:= \mathcal G^0[p^{\infty}]^{\f}_{\eta};\quad A[p^{\infty}]_{\eta}^{\tor}:=\mathcal A^0[p^{\infty}]^{\tor}_{\eta};\quad A[p^{\infty}]_{\eta}^{\f}:= \mathcal A^0[p^{\infty}]^{\f}_{\eta};$$ 
so that there are filtrations
$$G[p^{\infty}]_{\eta}^{\tor}\subseteq G[p^{\infty}]_{\eta}^{\f}\subseteq G_{\eta}[p^{\infty}]\quad \text{and}\quad A[p^{\infty}]_{\eta}^{\tor}\subseteq A[p^{\infty}]_{\eta}^{\f}\subseteq A_{\eta}[p^{\infty}].$$
\subsection{End of the proof}\label{conclusion}
By \cite[Expos\'e IX]{SGA7-II}, the inclusions $A[p^{\infty}]_{\eta}^{\tor}\subseteq A[p^{\infty}]_{\eta}^{\f}\subseteq A_{\eta}[p^{\infty}]$ produce a filtration of $A[p^{\infty}]_{\eta}$ giving semistable reduction for $A_{\eta}[p^{\infty}]$, in the sense that:
\begin{enumerate}
\item $A[p^{\infty}]_{\eta}^{\f}$ and $A_{\eta}[p^{\infty}]/A[p^{\infty}]^{\tor}_{\eta}$ extend to p-divisible groups $A[p^{\infty}]_{1}$ and $A[p^{\infty}]_{2}$ over $S$ (see \cite[Proposition  5.6, Pag. 380, Expos\'e IX]{SGA7-II}). 
\item If $f_A:A[p^{\infty}]_{1}\rightarrow A[p^{\infty}]_{2}$ denotes the natural induced map, then  $\Ker(f_A)$ is a multiplicative p-divisible group and $\Coker(f_A)$ is an \'etale p-divisible group (This follows from the orthogonality theorem \cite[Proposition  5.2, Pag. 372, Expos\'e IX]{SGA7-II}, which implies that $\Ker(f_A)\simeq\mathcal A^0[p^{\infty}]^{\tor}$ and $\Coker(f_A)\simeq (\mathcal (\mathcal A^{\vee})^0[p^{\infty}]^{\tor})^{\vee}$, where $\mathcal A^{\vee}$ is N\'eron-model of the dual abelian $A_\eta^\vee$ and $(\mathcal A^{\vee}[p^{\infty}]^{\tor})^{\vee}$ is the Cartier dual of $\mathcal A^{\vee}[p^{\infty}]^{\tor}$).
\end{enumerate}
To conclude the proof we now deduce for (1) and (2) above that the same properties holds for the filtration $G[p^{\infty}]_{\eta}^{\tor}\subseteq G[p^{\infty}]_{\eta}^{\f}\subseteq G_{\eta}[p^{\infty}]$.
\begin{enumerate}
\item By construction $G[p^{\infty}]_{\eta}^{\f}$ extends over $S$ to the p-divisible group $G[p^{\infty}]_1:=\mathcal G^0[p^{\infty}]^{\f}$. On the other hand, the diagram (\ref{diagramma}) shows that $G_{\eta}[p^{\infty}]/G[p^{\infty}]_{\eta}^{\tor}\simeq A_{\eta}[p^{\infty}]/A[p^{\infty}]_{\eta}^{\tor}$. So that we can set $G[p^{\infty}]_2:=A[p^{\infty}]_2$.
\item Let $f:G[p^{\infty}]_1\rightarrow G[p^{\infty}]_2$ be the induced morphism. We are left to prove that $\Ker(f)$ and $\Coker(f)$ are p-divisible groups and $\Ker(f)$ is multiplicative and $\Coker(f)$ is \'etale. This can be deduced from the analogues properties for $A[p^{\infty}]_1$ and $A[p^{\infty}]_2$, thanks to the commutative diagram with exact rows and columns:
\begin{center}
\begin{tikzcd}
0\arrow{r} & W[p^{\infty}]\arrow{r} \arrow[equal]{d} & \Ker(f)\arrow{r} \arrow{d}& \Ker(f_A)\arrow{r}\arrow{d} &0\\
0\arrow{r} & W[p^{\infty}]\arrow{r} & G[p^{\infty}]_1\arrow{r}\arrow{d} \arrow{d}& A[p^{\infty}]_1\arrow{r}\arrow{d} &0\\
 & & G[p^{\infty}]_2\arrow{r}{\simeq}\arrow{d} & A[p^{\infty}]_2\arrow{d}\\
 & & \Coker(f)\arrow{r}{\simeq}& \Coker(f_A).
\end{tikzcd}
\end{center}
\end{enumerate}
\bibliographystyle{alpha}
\bibliography{perfectfinalarxiv.bib}

\begin{thebibliography}{DdSMS91}

\bibitem[ABV05]{andreattabarbieri1motivi}
F.~Andreatta and L.~Barbieri-Viale.
\newblock Crystalline realizations of 1-motives.
\newblock {\em Math. Ann.}, 331:111--172, 2005.

\bibitem[AD22]{ambrosidaddezio}
E.~Ambrosi and M.~D'Addezio.
\newblock Maximal tori of monodromy groups of {$F$}-isocrystals and an
  application to abelian varieties.
\newblock {\em Algebr. Geom.}, 9(5):633--650, 2022.

\bibitem[BBM82]{berthelotbreenmessing2}
P.~Berthelot, L.~Breen, and W.~Messing.
\newblock {\em Th{\'e}orie de Dieudonn{\'e} cristalline {II}}, volume 930 of
  {\em Lecture Notes in Mathematics}.
\newblock Springer-Verlag, 1982.

\bibitem[Ber96]{berthelotrigide}
P.~Berthelot.
\newblock Cohomologie rigide et cohomologie rigide {\`a} supports propres
  (premi{\`e}re partie).
\newblock Preprint, 1996.

\bibitem[BL22]{supersingular}
D.~Bragg and M.~Lieblich.
\newblock Perfect points on genus one curves and consequences for supersingular
  {K}3 surfaces.
\newblock {\em Compos. Math.}, 158(5):1052--1083, 2022.

\bibitem[BLR90]{neronmodels}
S.~Bosch, W.~Lutkebohmert, and M.~Raynaud.
\newblock {\em N{\'e}ron Models}, volume~21 of {\em A Series of Modern Surveys
  in Mathematics}.
\newblock Springer-Verlag, 1990.

\bibitem[{\v C}es15]{kestutis}
K.~{\v C}esnavi\u{c}ius.
\newblock {P}oitou-{T}ate without restrictions on the order.
\newblock {\em Math. Res. Lett.}, 22:1621--1666, 2015.

\bibitem[D'A23]{daddeziominimal}
M.~D'Addezio.
\newblock Parabolicity conjecture of {$F$}-isocrystals.
\newblock To appear in Ann. of Math., 2023.

\bibitem[DdSMS91]{analyticprop}
J.~D. Dixon, M.~P.~F. du~Sautoy, A.~Mann, and D.~Segal.
\newblock {\em Analytic pro-p groups}, volume 157 of {\em London Mathematical
  Society Lecture Note Series}.
\newblock Cambridge University Press, 1991.

\bibitem[dJ95]{dejongcrystallineviaformal}
A.~J. de~{J}ong.
\newblock Crystalline {D}ieudonn{\'e} module theory via formal and rigid
  geometry.
\newblock {\em Publ. Math. Inst. Hautes \'Etudes Sci.}, 82:5--96, 1995.

\bibitem[dJ98]{dejonghomomorphism}
A.~J. de~Jong.
\newblock Homomorphisms of {B}arsotti-{T}ate groups and crystals in positive
  characteristic.
\newblock {\em Invent. Math.}, 134:301--333, 1998.

\bibitem[DK73]{SGA7-II}
P.~Deligne and N.~Katz.
\newblock {\em Groupes de monodromie en g\'eom\'etrie alg\'ebrique. {II}},
  volume 340 of {\em Lecture Notes in Mathematics}.
\newblock Springer-Verlag, 1973.
\newblock S{\'e}minaire de G{\'e}om{\'e}trie Alg{\'e}brique du Bois-Marie
  1967--1969 (SGA 7, {II}).

\bibitem[{\'E}te02]{etesseschemiab}
J.~Y. {\'E}tesse.
\newblock Descente {\'e}tale des {$F$}-isocristaux surconvergents et
  rationalit{\'e} des fonctions {$L$} de sch{\'e}mas ab{\'e}liens.
\newblock {\em Ann. Sci. Ec. Norm. Sup{\'e}r.}, 35(4):575--603, 2002.

\bibitem[Fal91]{Faltings}
G.~Faltings.
\newblock Diophantine approximation on abelian varieties.
\newblock {\em Ann. of Math.}, 133:549--576, 1991.

\bibitem[Ghi10]{ghiocaelliptic}
D.~Ghioca.
\newblock Elliptic curves over the perfect closure of a function field.
\newblock {\em Canadian Mathematical Bulletin}, 53(1):87--94, 2010.

\bibitem[GKU21]{cutbycurve}
T.~Grub, K.~Kedlaya, and J.~Upton.
\newblock A cut-by-curves criterion for overconvergent of {F}-isocrystals.
\newblock Preprint, 2021.

\bibitem[GM06]{ghiocamoosaperfection}
D.~Ghioca and R.~Moosa.
\newblock Division points on subvarieties of isotrivial semi-abelian varieties.
\newblock {\em Int. Math. Res. Not. IMRN}, 2006:1--23, 2006.

\bibitem[Hel22]{helmperfect}
D.~Helm.
\newblock An ordinary abelian variety with an {\'e}tale self-isogeny of p-power
  degree and no isotrivial factors.
\newblock {\em Math. Res. Lett.}, 20(2):445--454, 2022.

\bibitem[Hin88]{Hindry}
M.~Hindry.
\newblock Autour d'une conjecture de {S}erge {L}ang.
\newblock {\em Invent. Math.}, 94:575--603, 1988.

\bibitem[Hru96]{hrushovski}
E.~Hrushovski.
\newblock The {M}ordell-{L}ang conjecture for function fields.
\newblock {\em J. Amer. Math. Soc.}, 9(3):667--690, 1996.

\bibitem[Jan94]{jansen}
U.~Jannsen.
\newblock Motivic sheaves and filtrations on {C}how groups.
\newblock In {\em Motives ({S}eattle, {WA}, 1991)}, volume~55 of {\em Proc.
  Sympos. Pure Math.}, pages 245--302. 1994.

\bibitem[Ked04]{kedlayafully}
K.~S. Kedlaya.
\newblock Full faithfulness for overconvergent {$F$}-isocrystals.
\newblock In {\em Geometric Aspects of Dwork Theory}, volume~II, pages
  819--883. Adolphson et al, 2004.

\bibitem[Ked22]{kedlayanotes}
K.~S. Kedlaya.
\newblock Notes on isocrystals.
\newblock {\em J. Number Theory}, (237):353--394, 2022.

\bibitem[LN59]{langneron}
S.~Lang and A.~N{\'e}ron.
\newblock Rational points of abelian varieties over function fields.
\newblock {\em Amer. J. Math.}, 81(1):95--118, 1959.

\bibitem[Mes72]{messingcrystals}
W.~Messing.
\newblock {\em The Crystals Associated to {B}arsotti-{T}ate Groups: with
  Applications to Abelian Schemes}, volume 264 of {\em Lecture Notes in
  Mathematics}.
\newblock Springer-Verlag, 1972.

\bibitem[Mor19]{morrowvariational}
M.~Morrow.
\newblock A variational {T}ate conjecture in crystalline cohomology.
\newblock {\em J. Eur. Math. Soc.}, 21:3467--3511, 2019.

\bibitem[P{\'a}l22]{Pal}
A.~P{\'a}l.
\newblock The p-adic monodromy group of abelian varieties over global function
  fields of characteristic p.
\newblock {\em Documenta Mathematica}, 27:1509--1579, 2022.

\bibitem[R{\"o}s15]{rosslerperfectI}
D.~R{\"o}ssler.
\newblock On the group of purely inseparable points of an abelian variety
  defined over a function field of positive characteristic {I}.
\newblock {\em Comment. Math. Helv.}, 90:23--52, 2015.

\bibitem[R{\"o}s20]{rosslerperfectII}
D.~R{\"o}ssler.
\newblock On the group of purely inseparable points of an abelian variety
  defined over a function field of positive characteristic {II}.
\newblock {\em Algebra Number Theory}, 14:1123--1173, 2020.

\bibitem[Ser64]{serregroupesdecongruence}
J.P. Serre.
\newblock Sur les groupes de congruence des vari{\'e}t{\'e}s ab{\'e}liennes.
\newblock {\em Izv. Akad. Nauk SSSR Ser. Mat.}, 28:3--20, 1964.

\bibitem[{Sta}20]{stacks-project}
The {Stacks project authors}.
\newblock The stacks project.
\newblock \url{https://stacks.math.columbia.edu}, 2020.

\bibitem[Tat66]{Tatefinitefield}
J.~Tate.
\newblock Endomorphisms of abelian varieties over finite fields.
\newblock {\em Invent. Math.}, 2:134--144, 1966.

\bibitem[Tri08]{Trihan}
F.~Trihan.
\newblock A note on semistable {B}arsotti-{T}ate groups.
\newblock {\em J. Math. Sci. Univ. Tokyo}, 15:411--425, 2008.

\bibitem[Tsu23]{tsuzukiminimal}
N.~Tsuzuki.
\newblock Minimal slope conjecture of {$F$}-isocrystals.
\newblock {\em Invent. Math.}, 231:39--109, 2023.

\bibitem[Xin21]{yuan2021}
Y.~Xinyi.
\newblock Positivity of {H}odge bundles of abelian varieties over some function
  fields.
\newblock {\em Compos. Math.}, 157:1964--2000, 2021.

\end{thebibliography}

\end{document}